# Understanding the Role of Covariates in Numerical Reconstructions of Real-World Vehicle-to-Pedestrian Collisions

Natalia Lindgren[a], Svein Kleiven[a] and Xiaogai Li[a]

[a]*Division of Neuronic Engineering, KTH Royal Institute of Technology, Sweden,*

**ABSTRACT**

Traumatic Brain Injuries (TBIs) are a pressing global public health issue, impacting tens of millions of individuals annually. Vulnerable road users (VRUs), such as pedestrians, are vastly overrepresented in the worldwide TBI statistics. To evaluate the effectiveness of injury prevention measures, researchers often employ Finite Element (FE) models of the human body to virtually simulate the human response to impact in real-world road traffic accident scenarios. However, VRU accidents occur in a highly uncontrolled environment and, in consequence, there is a large amount of variables (covariates), e.g. the vehicle impact speed and VRU body posture, that together dictate the injurious outcome of the collision. At the same time, since FE analysis is a computationally heavy task, researchers often need to apply extensive simplifications to FE models when attempting to predict real-world VRU head trauma. To help researchers make informed decisions when conducting FE accident reconstructions, this literature review aims to create an overarching summary of covariates that have been reported influential in literature. The review provides researchers with an overview of variables proven to have an influence on head injury predictions. The material could potentially be useful as a basis for choosing parameters to include when performing sensitivity analyses of car-to-pedestrian impact simulations.

## 1. Introduction

Traumatic Brain Injuries (TBIs) are an urgent and world-wide public health concern, affecting up to 69 million individuals each year [1]. Such injuries may have fatal outcomes or severe effects on cognitive, physical and behavioral functions- returning to "normal baseline functioning" after such an injury can take months, if not years [2]. Global estimates show that out of the tens of millions of annual TBI occurrences, half can be linked to road traffic collisions [1]. Vulnerable road users (VRUs) - that is, pedestrians, cyclists, and powered two- and three-wheeler riders - are particularly overrepresented in the statistics [3]. In order to mitigate this public health burden, there is an evident need of targeted intervention strategies.

Evaluating the effectiveness of potential prevention measures is simply not possible without the fundamental knowledge on injury mechanisms. To characterize injury mechanisms under real-world conditions, researchers often turn to the indispensable procedure of accident reconstructions, in which laws of classical mechanics are used in combination with physical evidence to determine how and why an accident, and any sustained injuries, occurred. Employing computational models to reconstruct real-world accidents is a particularly popular option, as they offer a cost-effective, highly reproducible and robust alternative to the more traditional experimental approaches involving crash test dummies or human cadavers, i.e. Post Mortem Human Subjects (PMHS).

A widely used tool for computational accident reconstruction is Finite Element (FE) analysis, which is a numerical method routinely used in various engineering practices to predict the response of physical objects subjected to different types of loading. FE techniques can be used to create virtual, anatomically-detailed human surrogates, or so-called Human Body Models (HBMs), which in turn can be used to simulate the human response to impact. As opposed to dummies, an HBM or an FE head/brain model can offer measurements of tissue-based metrics for injury assessments, such as changes in strain of the brain tissue during loading (Figure 1). This enables researchers to get deeper insights into the injury mechanisms of specific tissues during dynamic loading scenarios, such as road traffic accidents [4].

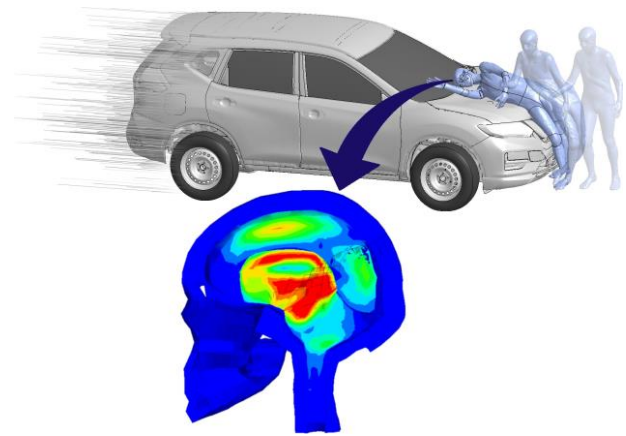

**Figure 1.** Reconstructions of real-world accidents using an HBM and/or FE head/brain model can be used to study the tissue-based mechanisms behind TBI, such as changes in strain of the brain tissue during loading.

Yet, using FE to reconstruct road traffic accidents is a challenging task. VRU collisions occur in a highly uncontrolled and unmonitored environment, and in consequence, there is a large amount of variables that together dictate the injurious outcome of the event. Some of these variables, or *covariates*, can heavily influence the dynamic response of the human body, the severity of injuries, and the overall body kinematics. Well-recognized influential covariates are for example the vehicle impact speed [5–7] and the VRU's posture prior to impact [8–10], although, judging by the current body of literature, the list of covariates is long.

*Corresponding author
E-mail address:* nalindgr@kth.se (N. Lindgren)

The reader should bear in mind that FE analysis is a computationally expensive method in which the computational cost grows with the level of detail of the model. Besides that, generating and modifying FE models is a notoriously time-demanding task. For these reasons, a modeler often strives to reduce the size and complexity of the problem, trying to balance the acceptable extent of accuracy with the time needed for pre-processing and computation. Hence, researchers are customarily bound to employ a wide variety of simplifications to FE models when investigating VRU collisions, models which are later used to draw conclusions about head injury mechanisms and injury risk. For example, some researchers disregard the influence of subject anthropometry and/or posture, and instead use HBMs with generic geometries [11–13], while others disregard the influence of cars' detailed design by using simplified or generic FE car models [13, 14].

At the same time, the output of an FE simulation is highly sensitive to boundary conditions and input parameters, and the reports on influential covariates in car-to-VRU collisions are piling up. The question arises whether the simplifications that are commonly adopted in the field are sufficient to capture the complex dynamics of a real-world head impact. Researchers may be making simplifications that are not motivated by the current state of knowledge in the field. To help researchers make informed decisions when carrying out computational accident reconstructions, an overarching summary of reported identified covariates is needed.

This report embodies a literature review aiming to summarize the current state of knowledge on covariates in car-to-pedestrian head traumas. The objective is to identify the variables reported to have an influence on head injury prediction in such impact scenarios, compiling which aspects that are of importance for performing accident reconstructions. Ultimately, this review aims to provide researchers with a list of covariates which would function as a general guideline, a checklist if you will, for performing VRU accident reconstructions using FE. The material could be used as a basis for choosing parameters to consider during accident reconstructions, as well as parameters to include when performing sensitivity analysis.

## 2. Review Method

This literature review is meant to be informative rather than all-compassing. A set of articles were chosen based on search words and search blocks, and then filtered out based on relevancy to the research question.

The literature search was primarily conducted using the KTH Library search service Primo. The following combination of keywords was considered as a base for the literature review: ("effect*" OR "influence*" OR "impact of") AND ("head" OR "brain") AND ("pedestrian*" OR "vulnerable" OR "VRU") AND ("car*" OR "vehicle*"). The search block resulted in 77 papers.

The research questions were answered by a targeted/focused literature review. Cited papers in key articles were also included in the review.

### 2.1. Delimitations

Pedestrian impacts can be divided into three phases. The first phase involves the initial contact phase, during which the pedestrian wraps around the front-end of the vehicle. The second phase is the flight phase, during which the pedestrian separates from the vehicle and is projected away of the vehicle. The third phase incorporates the rolling and sliding motion against the ground [15]. The third phase, involving secondary impacts to the ground, will be disregarded in this literature review. This is firstly because of the potential computational hindrance this would introduce in an FE reconstruction- FE simulations of the ground-strike would demand significantly longer simulation time. Secondly, several researchers have pointed out how seemingly unpredictable post-impact kinematics are, even at low speeds, due to the complexity of the many variables involved in a crash [5, 16]. Thirdly, many researchers have reported that ground contact account for a small percentage of head injuries in car impact cases [6, 17, 18]. Reported head injuries seem more likely to be sustained from impacts with structural parts of the impacting vehicle than the ground [18].

Of course, the modeling of the human body and head/brain, including the constitutive models of tissues, part constraints, the mesh, contact settings and so forth, has an important influence in the prediction of head injuries. However, detailed head/body FE modeling is excluded from the scope of this study.

## 3. Terminology

In this section, short explanations on the vocabulary used in the processed articles will be provided.

### 3.1. Injury metrics

In the reviewed articles, the risk of head injury is often evaluated in terms of a chosen quantitative risk metric, or an injury criterion. An injury criterion correlates a physical parameter, such as acceleration, with the probability of head injury. In automotive safety, head acceleration measurements using crash test dummies have conventionally been used in head injury assessment [19]. The most widely used and acknowledged criteria is the so-called Head Injury Criterion (HIC), which was originally based on the assumption that the head's Center of Gravity (CoG) linear acceleration and its duration alone is an indicator of injury

[20]:

$$HIC \max_{t_1,t_2}\left[(t_2 - t_1)\left(\frac{1}{t_2 - t_1}\int_{T-1}^{t_2} a(t)\,dt\right)^{2.5}\right]$$

Here, $t_2 - t_1$ is the critical time period of deceleration during impact (usually 15 ms, also known as HIC15) and $a(t)$ is the resultant deceleration of the head CoG at time $t$. For reference, if HIC exceeds 1000 the case is commonly regarded as life threatening [21].

The brain has for long been acknowledged to be sensitive to rotations [22]. Since HIC does not take rotational motion into account, authors have presented a variety of different global injury criteria. The Generalized Model for Brain Injury Threshold (GAMBIT) [23], the Head Impact Power (HIP) [24], the Brain Injury Criterion (BrIC) [25], and the Universal Brain Injury Criterion (UBrIC) [26], are some examples of suggested kinematics-based injury criteria, many of which include terms of rotational kinematics.

There are also many tissue-based injury metrics for assessing brain injury risk. Examples of such criteria include the Maximum Principal Strain (MPS), which measures the maximum tensile strain experienced by the brain tissue during an impact, and the Cumulative Strain Damage Criterion (CSDM), which is a measurement of the accumulative volume of brain tissue that endures a specific level of strain.

A common metric used for quantifying the translational movements of the head during impact is the head's Peak Linear Acceleration (PLA), while Peak Angular Velocity (PAV) and Peak Angular Acceleration (PAA) are used when assessing rotational head kinematics. Within this context, PLA, PAV and PAA are always measured relative to the head's CoG.

### 3.2. Common assessment methods

The reviewed studies often choose to evaluate standardized tests to simulate the most common pedestrian-vehicle crashes. Test specifications and rating systems for assessing the pedestrian injury potential of vehicle front structures has previously been developed by the European Enhanced Vehicle-Safety Committee (EEVC), and are also included in the European New Car Assessment Program (EuroNCAP) car rating program. The ratings involve four subsystem impactor tests: adult headform impacts against the hood or windshield base, smaller headform impacts representing children against the hood, and upper legform impacts against the Bonnet Leading Edge (BLE) and bumper, see Figure 2. These impact configurations are frequently employed as they provide standardized, repeatable measures of vehicle safety in real-world pedestrian accidents.

The tests are designed to mimic a 40 km/h car-pedestrian impact with a pedestrian moving laterally across the path of the car. 40 km/h is chosen since it is considered to be a critical speed for the onset of serious/fatal injuries [27], and is thus a frequently evaluated impact velocity.

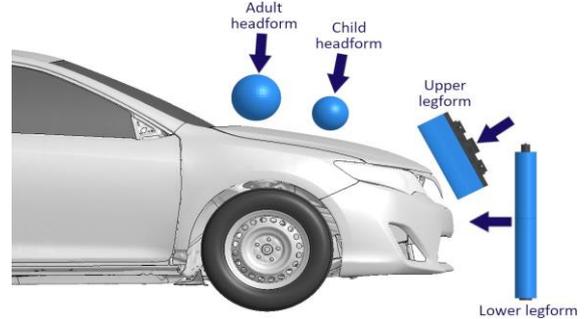

**Figure 2.** EuroNCAP pedestrian safety test protocol.

Apart from HBMs and FE head models, researchers often employ Rigid Body Models (RBMs) to simulate pedestrian accidents. RBMs consist of rigid (undeformable) bodies joined together to mimic the overall human structure, mass, mass distribution, and joint movements. Unlike HBMs, RBMs cannot model tissue deformation. Yet, they are well-suited for parametric studies, as they can provide relatively good accuracy at a low computational expense.

Researchers often use CORrelation and Analysis (CORA) scores to compare time-history signals, such as head impact velocity over the duration of an impact, between test data and simulations. A CORA score of 0 indicates no correlation of the pulses, while 1 indicates a perfect near-perfect correlation. The CORA rating system is particularly common when comparing FE models to PMHS experiments.

### 3.3. Car terminology

A car involves many structural components, of which the front-end structure of the vehicle is of most significance for pedestrian impacts. In Figure 3, important car features are named and illustrated. The Wrap Around Distance (WAD) and Bonnet Leading Edge (BLE), which are consistently referred to in this report, are illustrated in the figure as well.

Passenger vehicles can be divided into several classifications based on their use and characteristics. The geometry, in particular the front-end profile of vehicles, varies among car classes. They usually have distinctly different BLE heights, BLE shapes, bonnet lengths and hood/windshield inclinations. Representative front-profiles of a sedan, Sports Utility Vehicle (SUV) and a van are shown in Figure 3.

The *yaw angle* describes the orientation of the car relative to a fixed reference direction, usually the direction of the road or a global reference frame (e.g., north). It represents the car's rotational position around its vertical (z) axis. A change in yaw angle indicates that the car is rotating or "yawing" to the left or right. The *steering angle* is the input provided by the driver (via the steering wheel) to control the direction of the front wheels.

*Pitching* of a vehicle occurs during braking and refers to a rotational motion of a car's body which causes the front end to dip downward and the rear end to rise. Pitching occurs due to braking deceleration, transferring the car's weight from the rear to the front wheels.

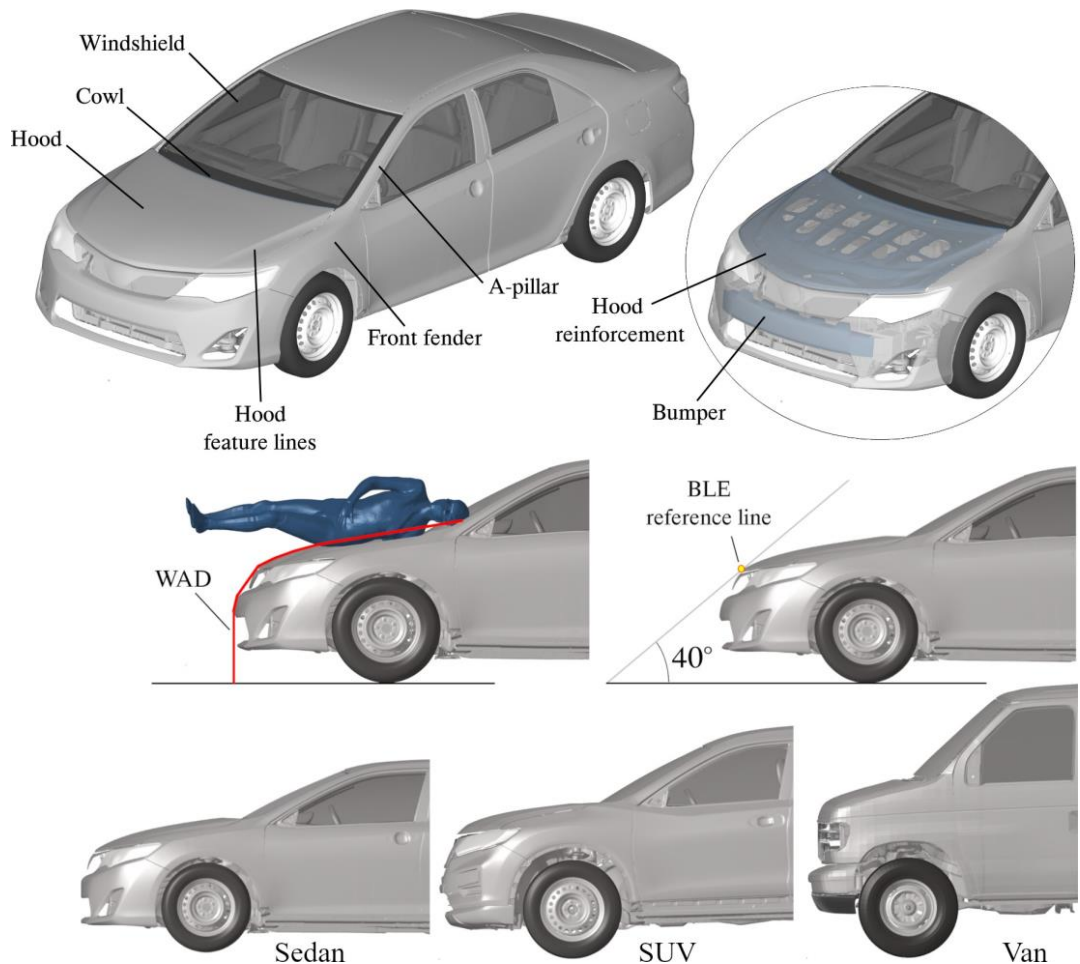

**Figure 3.** Terminology used in the reviewed literature.

The *hood hinge* is typically located at the rear corners, below the hood, and are the structural components of the car that allows the hood to open and close. The *hood stoppers* are the cushioned resting points for the hood when closed. The *cowl* refers to the structural component located below the windshield base and the hood edge of a vehicle. This area typically houses heating and ventilation systems, and contains the drainage channels and windshield wipers and wiper arm linkages.

## 4. Results: Identified covariates

More than 40 different covariates, reported to influence head impact kinematics and injury outcome in car-to-pedestrian collisions, have been identified in this review, see Figure 4. The covariates were broadly divided into six main categories. The first category of covariates relates to the *front-end geometry of the impacting vehicle*. The second category, *local car stiffness*, relates to the stiffness of the impacting surface of the car, at the region of head impact, while the third category, *global vehicle stiffness*, relates to vehicle stiffness not necessarily at the region of head impact. The fourth category relates to *subject anthropometry*, while the fifth relates to the *subject's pre-impact conditions*. The final categories relate to the *impact kinematics* and applied *boundary conditions*.

Most of the listed covariates have a considerable influence on the head impact location. The impact location has, judging by the body of literature processed in this review, a large influence on the likelihood of injury. Many of the covariates also influence the kinematic motion of the head, which is important to consider in car-to-pedestrian impacts, since it not only determines the head impact location, but also the head impact velocity, rotational kinematics and forces behind the head strike.

In the following subsections, the listed covariates will be explained in more detail.

### 4.1. Vehicle front-end geometry

There seems to be a consensus in the field that one of the main factors dictating the head impact characteristics is the geometry of the vehicle front-end. The vehicle geometry dictates where on the vehicle surface the head impacts, which in turn, as further explained in Section 4.2.1, can determine the forces and kinematics sustained by the pedestrian's head. There are many factors that together make up for the vehicle front-end geometry, and some of these factors are listed in this section.

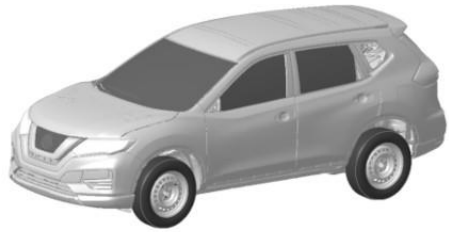
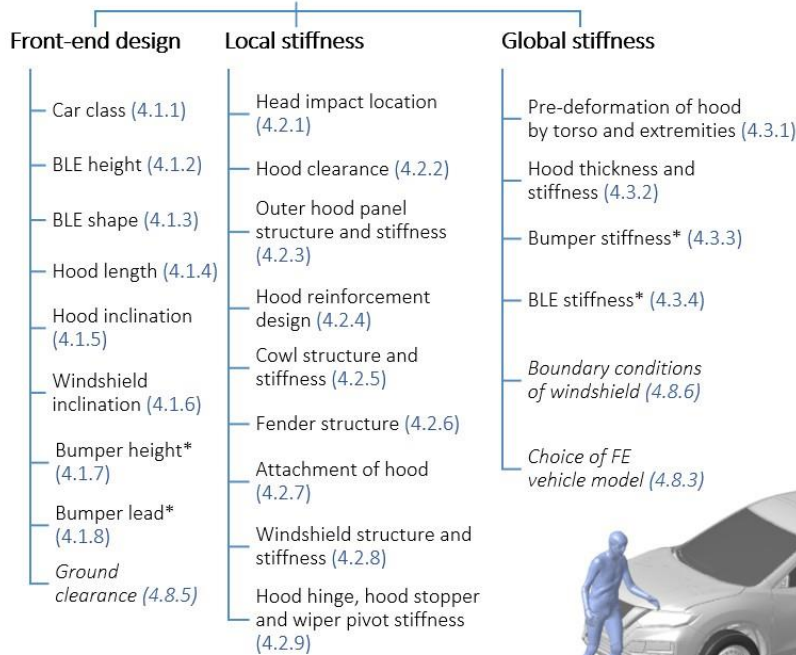
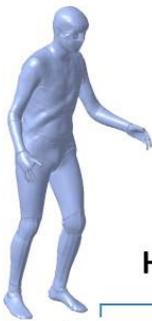
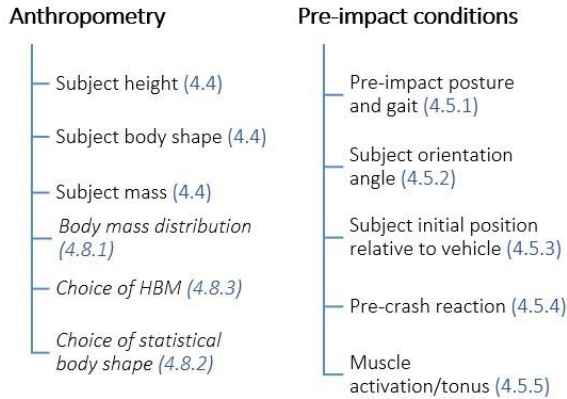
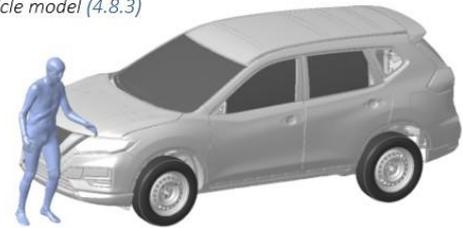
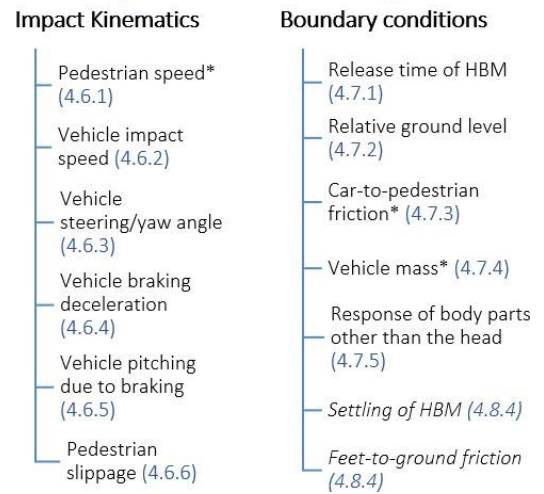

**Figure 4.** Influential variables affecting head impact response in car-to-pedestrian impacts, both in terms of head kinematic response and head injury prediction. The factors written in cursive have been identified by the authors as potential covariates, but without any findings in literature. They can be interpreted as knowledge gaps within the field.
*Previous publications suggests that these covariates are not particularly influential for head impact prediction.

It is generally believed that the shape of the vehicle front-end influences the energy associated with the car-to-pedestrian collision, while the general stiffness (scrutinized in Section 4.2 and 4.3) of the vehicle might determine the corresponding force of the impact [28].

*4.1.1. Car class*

The main features of the front-end geometry will be determined by the classification of the vehicle. There is a convincing amount of evidence suggesting that the vehicle class, and conversely the front-end profile geometry, of the impacting vehicle in a car-to-pedestrian impact has an immense effect on head injury outcome [6, 7, 18, 28–47]. Mainly, the vehicle class will influence where on the vehicle the head strikes: due to the differing BLE heights, SUVs and light trucks pose a greater risk of head impacts to the hood, while sedans and other passenger vehicles pose a greater risk of windshield head impacts [18, 43]. This difference can be depicted by the experimental findings of Kerrigan et al. [48], who performed full-scale PMHS impacts comparing lateral impacts with a small sedan to impacts with a large SUV. High-speed video frames of the head impacts (Figure 5) clearly shows how the highly situated hood of an SUV results in head-to-hood contact rather than head-to-windshield.

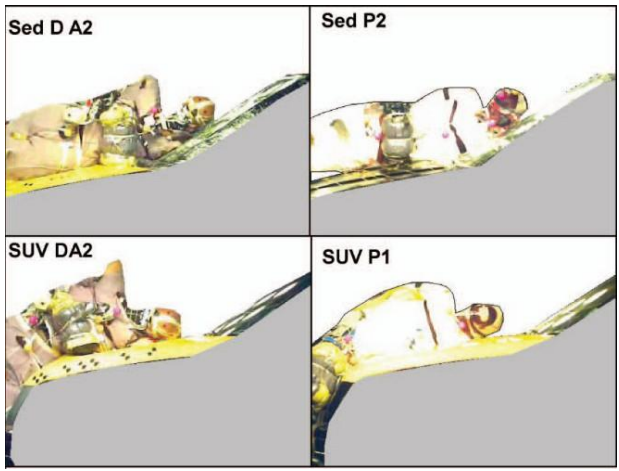

**Figure 5.** High-speed video snapshots of head impacts of a crash dummy (left column) and PMHS (right column) impacted by a small sedan (top row) and a large SUV (bottom row). Figure adapted from Kerrigan et al. [48].

There is evidence suggesting that some car classes also pose a greater hazard to pedestrians than other cars [6, 18, 38, 48]. For instance, SUVs have been shown to increase injury and fatality risk among pedestrians compared to smaller, passenger vehicles [6, 38]. Researchers generally attribute this to the difference in BLE height among car classes, which is an important covariate of head injury likelihood, see further details under Section 4.1.2.

There might be stiffness differences between different car classes, particularly since stiffness is generally believed to correlate with mass [6]. The differing stiffness properties of vehicle classes can be illustrated by the findings of Martinez et al. [49], who summarized the results of 69 sets of physical legform and headform impact tests against vehicle fronts in form of average stiffness mappings of different car classes, see Figure 6. SUVs resulted in notably higher headform HIC scores from hood and windshield impacts compared to other vehicle classes, suggesting a stiffness variation among vehicle classes.

Simms et al. [6] claim that the windshield construction of SUVs compared to sedans are similar.

*4.1.2. BLE height*

The influence of the height of the BLE on head injury has been extensively researched and is undoubtedly, judging by the many published studies on the topic, one of the main vehicle features influencing the head impact characteristics in car-to-pedestrian collisions [6, 28, 29, 36, 45, 50–54]. In particular, the BLE height, which is strongly related to the vehicle class, seems to dictate where on the vehicle the head strikes. For impacts with adult pedestrians, a car with a high BLE and long hood result in head contact on the bonnet top, while a short BLE and a short hood tend to result in windshield impacts [43]. Seemingly, the BLE height affects how the body rotates about the pelvis and how much of the pedestrian's body mass is engaged in the initial impact. This influences the energy being transferred to the head during collision.

For example, Shi et al. [52] performed 252 car-to-pedestrian impact simulations using RBMs, including a range of vehicle front-end designs and impact velocities, and could demonstrate that the BLE height was the leading governing factor of the head impact location. Moreover, Anderson and Doecke [29], showed that a high BLE increased neck loads in car-to-pedestrian collisions, which could lead to amplified HIC estimations. The study was based a large set of RBM car-to-pedestrian simulations (324 in total), involving a variation of vehicle geometries, pedestrian orientations, impact speeds and braking levels.

*4.1.3. BLE shape*

Some researchers have pointed out that the shape of the BLE has an influence on head PLA and PAA in car-to-pedestrian collisions. For instance, Tolea et al. [45] used RBMs to investigate how changes in geometric parameters of the vehicle front-end, such as the BLE height, hood length, hood inclinations and BLE radius, affected the head kinematics of an impacting pedestrian. A significant influence of the BLE radius, describing the roundness/sharpness of the BLE, on the head impact PAA and PLA was observed. A set of different car classes was included in the analyses, including sedans, SUVs and vans.

Furthermore, in lateral vehicle impacts of PMHS, several authors have observed the occurrence of sliding/slipping of the subject over the bonnet during the contact phase. Kerrigan et al. [55] hypothesizes that such sliding up the hood by PMHS is promoted by a smooth sloping shape of the hood. It has been acknowledged that sliding has an influence on head impact kinematics [36, 56–58] (details on sliding is included in Section 4.6.6).

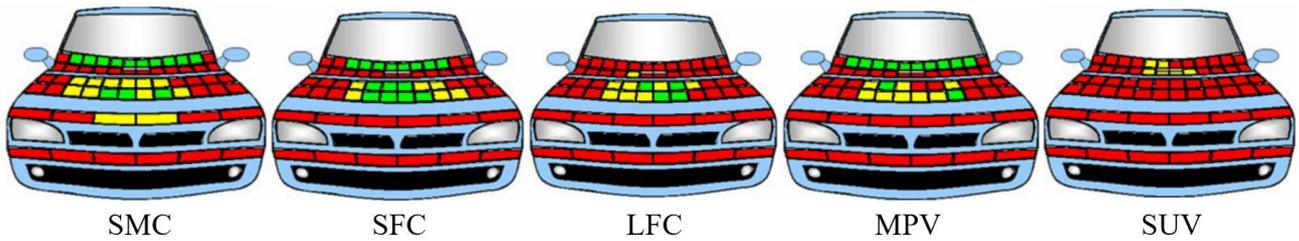

**Figure 6.** Mappings of average headform HIC scores and maximum legform bending of different vehicle segments for a variation of car classes: super minicars (SMC), small family cars (SFC), large family cars (LFC), multipurpose vehicles (MPV) and sport utility vehicles (SUV) (left to right). Red areas on the hood and windshield indicate high headform HIC score (<1350), while green indicate low HIC scores (<1000). Figure adapted from Martinez et al. [49] with permission from the authors.

### 4.1.4. Hood length

Together with the BLE height, the hood length can dictate whether or not the head impact will be on the hood or on the windshield [28, 43, 45]. Generally speaking, a passenger vehicle with a long hood would result in a head impact on the hood, in the same way as a vehicle with a highly situated BLE would. Correspondingly, a vehicle with a short hood would likely result in a head impact against the windshield [18, 28, 59].

### 4.1.5. Hood inclination

Ahmed et al. [60] simulated headform impacts on three points on a vehicle hood with three different hood inclination angles (6°, 8°, 10°). For a head impact on the centerline and middle of the hood, the estimated HIC varied significantly, where an inclination angle of 10° resulted in a HIC score of 1900 and an inclination angle of 6° resulted in a HIC of 1200, indicating that the hood incline might be very influential in the response of head-to-hood impacts. A number of other researchers have brought up the hood inclination as an influential parameter for predicting head impact response [54, 61].

### 4.1.6. Windshield inclination

Several researchers have identified the inclination angle of the windshield to be influential to head impact kinematics [62–64]. For instance, Wang et al. [64] drew this conclusion based on pedestrian FE-coupled RBM simulations, simulating 45 km/h lateral head-to-windshield impacts. The authors reported that both HIC and PAA were significantly affected by the inclination of the impacting windshield, which varied between 24° and 50° (evaluated with 2° intervals). Within the evaluated interval of angles, the HIC varied between 360 and 1200, while the PAA ranged between 10,000 and 23,000 rad/s$^2$. Similar results were observed by Lyons et al. [65], who studied the influence of the windshield angle (20° to 55°) in sedan-to-pedestrian impacts at 40 km/h using RBMs. It was shown how an increase in windshield angle could reduce the head PAA by up to 18%, see Figure 7.

### 4.1.7. Bumper height

Yang et al. [53] performed a sensitivity study involving lateral car-to-pedestrian collisions at different impact velocities and different bumper heights, simulated using an pedestrian FE-coupled RBM. The researchers were able to show that the estimated HIC was almost unaffected by the bumper height (ranging between 275 to 500 mm above ground level). These findings are supported by other studies, see e.g. [66] or [36].

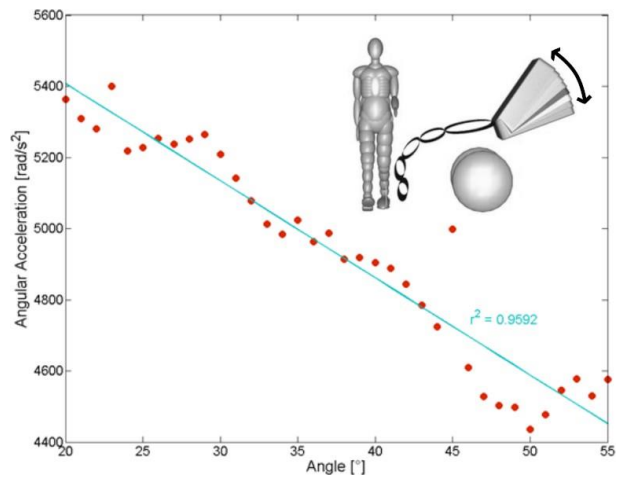

**Figure 7.** Variation of PAA with windshield inclination in 40 km/h lateral sedan-to-pedestrian impacts using RBMs. Figure adapted from Lyons et al. [65] with permission from the authors.

### 4.1.8. Bumper lead

In a parametric study, Li et al. [36] investigated the influence of various shape parameters on the prediction of head injuries in lateral impacts at an impact speed of 40 km/h using RBMs. The researchers reported no considerable effect of the bumper lead (the protrusion of the bumper) on the head impact velocity in simulations. Bumper leads of 50, 100 and 200 mm were included in the analysis.

## 4.2. Local vehicle stiffness

This section relates to covariates that are associated with the stiffness of the impacting surface of the car at the region of head impact.

### 4.2.1. Impact location

The impact location of the head is a major determinant of head injury severity, mainly due to the variation in stiffness properties over different locations on the vehicle front-end and its underlying structures [32, 36, 39, 40, 43, 45, 67, 68]. Many authors point out some regions of the car that are notably stiffer than others, and thus pose a greater risk of head injury. These regions are the windshield edge and bottom corners [6, 28, 39, 43, 45], A-pillars [6, 28, 39, 43], cowl [6, 28, 39, 43], hood/fender seam [28, 39, 69], hood hinge and stopper [28, 39, 43], and wiper pivot [39]. Some of the more compliant parts of the vehicle, which also pose less head injury risk, are the center of the windshield and, if the under-hood clearance is large enough, the center of the hood [39].

To illustrate the vehicle stiffness variations of a vehicle front-end, the results from a set of standard headform experiments carried out by Mizuno and Kajzer are provided [70], see Figure 8. The headform tests were aimed to replicate a 40 km/h pedestrian impact, and the impact angle was set to 65°, following standardized pedestrian test procedures.

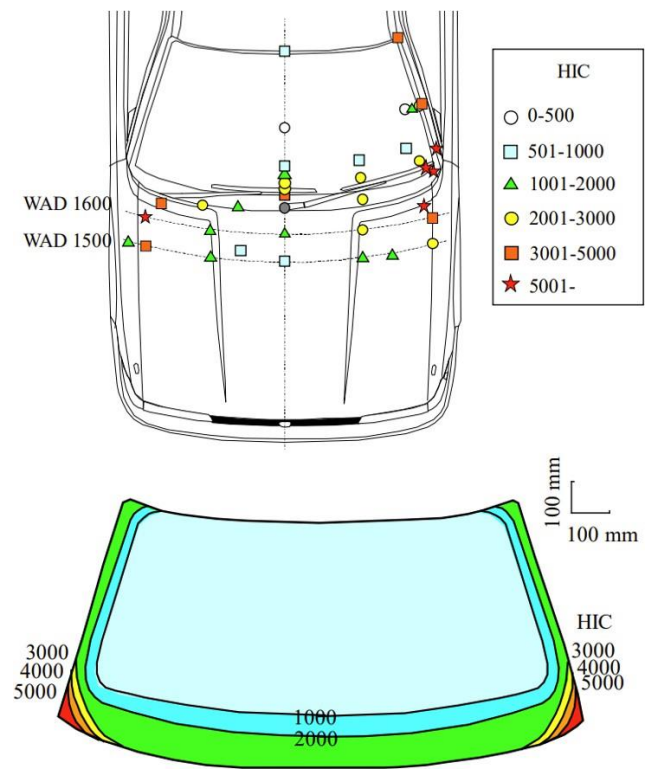

**Figure 8.** Variation in HIC scores of Japanese cars for various head impact locations, based on headform impacts representing 40 km/h impacts performed. Figure borrowed from Mizuno and Kajzer [39].

In line with the presented work of Mizuno and Kajzer [70], Watanabe et al. [68], who performed 72 FE simulations using a variation of vehicle models, pedestrian models (5th, 50th and 95th percentile) and vehicle speeds, observed that the head contact force was strictly related to the impact location on the vehicle. High impact forces were estimated when the HBM head impacted metal parts such as the hood and A-pillar, compared to when the head impacted the windshield glass. It was concluded that different vehicle regions pose different risk of head injury. Supporting findings were presented by Longhitano et al. [37], who analyzed pedestrian crash data covering pedestrian crashes in six US cities between the years of 1994 and 1998. The researchers were able to link the pedestrian risk of severe injury to different regions of the car, see Figure 9.

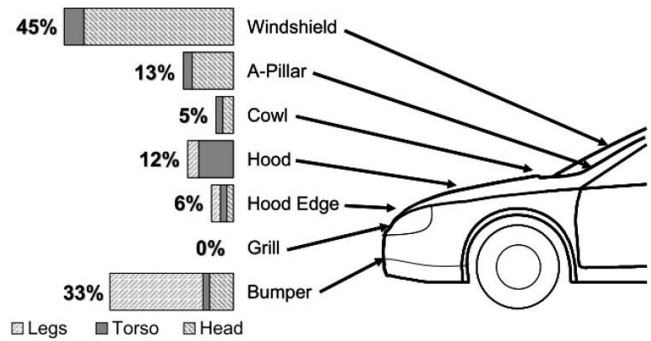

**Figure 9.** Distribution of AIS 3+ injuries for passenger cars. AIS 3+ refers to injuries classified as serious, severe, critical, or unsurvivable. A majority of severe head injuries were due to contact with the windshield. Figure borrowed from Longhitano et al. [37].

### 4.2.2. Hood clearance

Studies have shown how the hood clearance, which is the distance between the hood and the stiff underlying engine components, has a large influence on the head injury likelihood in head-to-hood impacts [43, 71–76]. The available hood deformation space can hinder a pedestrian from contacting underlying components, which are usually very stiff and, thusly, particularly hazardous. Fredriksson et al. [71] studied free falling headform impacts on a physical hood model with an underlying, adjustable plate mounted underneath. It was demonstrated how an increase of under-hood distance from 20 mm to 100 mm could reduce the HIC scores from 5000 to 500. The increase in under-hood distance was also shown to reduce the estimated CSDM from 20% to 0%. In conclusion, contact with stiff engine structures has a significant contribution to head injuries. This occurs if there is an insufficient clearance between the outer hood at the region of head impact.

### 4.2.3. Outer hood structure and stiffness

Several researchers have emphasized how the structure and stiffness of the hood is influential for HIC predictions in head-to-hood impacts [43, 77–79]. The hood design varies among car manufacturers and car models, and design choices such as feature lines can have an effect on the predicted head impact response.

Using full-body HBMs to simulate head-to-hood impacts, Chen et al. [77] demonstrated how stiffness changes (within reported ranges) in terms of hood thickness, elastic modulus and yield stress of the outer hood panel had an influence in estimated HIC values. On average, the HIC was more than doubled due to a thickness increase in both the inner and outer hood panel of around 1 mm. However, in most scenarios involving a hood stiffness increase or decrease, no statistically significant change in BrIC was observed.

By simulating headform impact response in hood impacts, Nie et al. [80] showed how head accelerations, both in terms of HIC and PLA, can increase significantly if the head impacts a hood feature line. The feature lines change the available deformation space for the head, influence the head's rotation and rebound.

### 4.2.4. Hood reinforcement design

As for the outer hood layer, the inner hood layer varies in design among car models. An example of variations in designs of the inner hood panel, or the hood reinforcement, is presented in Figure 10. Several researchers point out that the design of the hood reinforcement has an influence on head kinematics [67, 69, 74, 79]. It is obvious that head impacts close to a "rib" of a hood reinforcement differs in head impact response compared to head impacts close to a "hole" in the reinforcement. Similar to the hood clearance (see Section 4.2.2), this can be attributed the available crush space. When the crush space is small, e.g. at a hood reinforcement rib, the head impact will result in higher HIC scores compared to when the crush space is large, e.g. at a hood reinforcement hole.

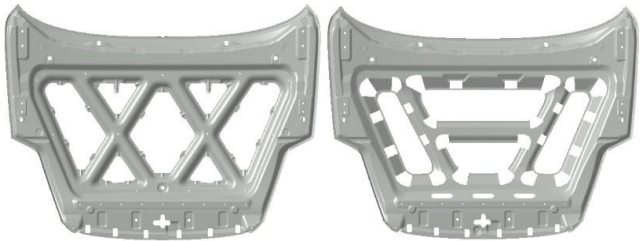

**Figure 10.** Examples of hood reinforcement structure designs. Figure borrowed from Lv et al. [69].

### 4.2.5. Cowl structure and stiffness

By FE analysis, researchers have found that head impacts on the cowl top of the vehicle result in significantly higher HIC values compared to other parts of the vehicle [69, 81], and that the stiffness of the cowl has an influence on head impact response in those impacts [82]. Other researcher stress that the hood overhang over the cowl might be of importance in head impacts near the cowl region [43]. A portion of severe head injuries have been attributed to the vehicle cowl (see Figure 9), thus this region of the car should be modeled with the same care as the windshield and hood.

### 4.2.6. Fender structure

In one study, the fender structure has been deemed influential [69]. The parting, where the hood continues to the fender, is reported to be notably stiff, which is also visualized by Mizuno and Kajzer [70] (recall Figure 8). According to Lv et al. [69], the stiffness of the fender parting could be explained by the reduced crush space close to the fender brackets.

### 4.2.7. Attachment of hood

Otubushin et al. [67] points out that the constraint between the inner hood panel and the outer hood panel is influential in head injury outcome. By numerical simulation, the researchers showed that having no constraints versus connecting the two hood layers with rigid links changed peak accelerations with 10% in headform impacts to the hood. The boundary support of the hood has also been shown to affect HIC values in full-scale lateral HBM head-to-hood impacts, see Section 4.2.9.

### 4.2.8. Windshield structure and stiffness

Chen et al. [77] compared two different windshield stiffnesses in lateral car-to-pedestrian impacts. In one configuration, the yield stress of the outer glass layers and the elastic modulus of the plastic interlayer was set to 20 MPa respective 1.4 GPa, while in the other configuration the corresponding values were set to 4 MPa and 1 GPa. The researchers concluded that softening and stiffening of the windshield had a statistically significant effect in the predicted HIC value (8% difference in HIC), whilst the BrIC was left essentially unaffected (BrIC values were estimated to 1.386 and 1.384).

Munsch et al. [83] implies that the design factor related to the windshield which has the most effect on HIC is the glass thickness. The thicker the glass, the higher the magnitude of HIC. Munsch et al. [83] suggest that the yield stress of the glass and the PVB, and the PVB thickness has negligible influence on PLA and HIC in comparison.

How the windshield is modeled in the FE simulation has shon to be very important. Alvarez et al. [84] compared two modeling approaches of windshields in head-to-windscreen impacts. One modeling approach involved modeling the two layers of glass and intermediate layer of PVB as two coinciding shell layers, of which one was able fracture. The second technique involved a three-layered-shell modeling approach, separated by the distance of thicknesses. The authors emphasize how the choice of windshield modeling had a large effect on the predicted brain strains.

Munsch et al. [83] further showed that the constitutive modeling of the windshield might be important when predicting head injuries in pedestrian collisions. The researchers conducted physical experiments with headforms impacting a vehicle side window, and demonstrated that three-layered laminated glazing, composed by two glass-layers separated by a plastic interlayer, produces significantly lower head injury risk compared to tempered glass (bare in mind, that tempered glass has gradually been replaced by laminated glass, and is not really used in modern vehicles).

Simms et al. [6] reported, based on a set of RBM simulations of frontal and lateral pedestrian impacts at different impact speeds and using different car classes, that the pedestrian head loading in head-to-windshield impacts is governed by the deformation behavior of the windscreen. As the windshield breaks, the stretching of the plastic laminate acts to substantially increase the estimated HIC scores. The forces affecting the head will be determined by how and when the glass fails, and how much the laminate is stretched.

### 4.2.9. Hood hinge, hood stopper and wiper pivot stiffness

Mizuno and Kajzer [39] performed headform impacts to a car, simulating a 40 km/h pedestrian collision, and found that the locally high stiffness of the hood hinge, hood stopper and wiper pivot had a major effect on HIC. The HIC was extremely high (>5000) for impacts on the hood hinge and the hood stoppers, see Figure 8. A total of 38 impact tests were performed on different areas of the car, including the cowl.

Chen et al. [77] investigated how the head impact response was influenced by strengthening or weakening of the hood boundary support in full-scale lateral HBM head-to-hood impacts. The boundary support was altered by modifying the thickness and material parameters of hinges on both sides of the hood. It was found that an increased stiffness of the hood hinge supports had a statistically significant effect on HIC scores, although the head did not directly impact the hood hinge in these simulations.

## 4.3. Global vehicle stiffness

This section relates to covariates that are associated with the stiffness of the impacting surface of the car, not necessarily at the region of head impact.

### 4.3.1. Pre-deformation of hood by torso and extremities

Chen et al. [77] performed 54 full-scale HBM simulations of lateral car-to-pedestrian impacts, including variations in HBM posture, orientation, anthropometry, impact location and vehicle stiffness in the test matrix, and emphasized that the pre-deformation of the hood caused by the torso, shoulder and upper extremities caused changes in head response that was not reflected in head-only impacts using a headform impactor.

### 4.3.2. Outer hood thickness and stiffness

In Section 4.2.3, it was shown how head impacts to the hood of a vehicle is highly influenced by the compliance of the hood. But the stiffness properties of the hood has also shown to influence head response in head-to-windshield impacts. Simms et al. [6] used RBMs to simulate a lateral head-to-windshield impact at 54 km/h, demonstrating how a 20% stiffer hood increased the estimated HIC by 20% (from 2200 to 2600).

### 4.3.3. Bumper stiffness

Liu et al. [36] investigated the influence of various shape parameters on the prediction of head injuries in lateral impacts at an impact speed of 40 km/h. The authors reported that the bumper stiffness affected the WAD by delaying the time of head contact. However, the head impact velocity was not notably affected.

### 4.3.4. BLE stiffness

In studies performed using RBMs, the stiffness of the BLE has been shown to have little effect on the head kinematics [9, 36, 50]. For example, Elliott et al. [50] varied the BLE stiffness in simulations of lateral car-to-pedestrian impacts and found that the head impact was largely unaffected by doubling or halving the BLE stiffness (1-2% difference in WAD and head rotation). Other authors confirm that there have been no reported relationships between the stiffness of the BLE and HIC scores [28, 43].

## 4.4. Pedestrian anthropometry

In recent decades, a large amount of research has been published, reporting on the importance of accounting for subject anthropometry for injury prediction [46, 55, 58, 68, 85–89].

Yan et al. [89] showed the importance of accounting for subject anthropometry by comparing the impact response between two generic HBMs, representative of two different populations. The statistically-derived body and skeleton morphology of the HBMs represented a 50$^{th}$ percentile American/European male HBM (THUMS AM50) and a 50$^{th}$ percentile Chinese male. In lateral impacts with different car classes, the authors could report significant differences in the predicted severity of brain and skull damage between the two models. For instance, in identical lateral 40 km/h impacts with a sedan, the predicted MPS between the who models were 23% higher for the Chinese HBM (0.86 versus 0.7), suggesting that pedestrian anthropometry, in particular among different populations and ethnicities, should not be neglected. The two models differed 9.6 cm in height and 10 kg in weight, and the authors do not report on whether the most influential factor was the height/size, weight or body morphology of the HBM.

Furthermore, it has been suggested that subject height has a particularly large influence on the head impact location and the head impact velocity [40, 42, 51, 58]. Schroeder et al. [42], who looked at PMHS tests, observed that shorter subject's heads tended to collide with the hood, while the taller subjects more often collided with the windshield. Subit et al. [88] compared the impact kinematics in 40 km/h lateral impacts of four full-scale PMHS of different statures. Two PMHS were targeted to be relatively short subjects (154 and 161 cm), while the other two were relatively tall (182 and 183 cm). It was shown that HIC and PLA increased with body height. Short PMHS produced shorter WAD, while tall PMHS experienced more sliding.

The fact that taller pedestrians experience head contact farther up the vehicle, as well as more sliding, has also been demonstrated in previous PMHS tests by Kerrigan et al. [55]. Kerrigan et al. suggested that the stature of the pedestrian could directly result in a difference in head injury risk, as a taller pedestrian risks contact with stiffer structures of the vehicle, such as the cowl. The findings also suggest that the mass of the pedestrian has an influence on hood deformation at pelvis contact.

Paas et al. [86] compared different scaling techniques of HBMs and concluded that adding a scaling factor to adjust the HBM stature and mass improved the prediction of the WAD, head impact velocity, trajectories and overall CORA rating, in comparison to a generic, unscaled HBM. Per example, for one PMHS experiment, the CORA rating was increased from 0.50 to 0.63 by global scaling of the HBM to match the PMHS stature. Kerrigan et al. [55] has previously shown how a global scaling factor of HBMs is not sufficient to predict pedestrian PMHS response in lateral impacts, indicating that Body Mass Index (BMI, calculated by dividing the subject's mass with their squared height) and overall body shape has an influence, and not solely the body height and mass.

Chen et al. [85] reconstructed PMHS impacts of lateral impact tests using two PMHS targeting an obese population. By comparing geometrically morphed HBMs with a generic HBM, the authors could demonstrate how poorly the generic-shaped HBMs captured the PMHS head kinematics and injury risks, whereas the morphed HBMs reproduced biofidelic responses in terms of head trajectories, velocities, accelerations and strains.

## 4.5. Pedestrian pre-impact conditions
### 4.5.1. Pre-impact posture and gait

The pedestrian posture has been shown to be highly influential for head impact kinematics in pedestrian collisions [8–10, 41, 46, 50, 51, 58, 90]. For instance, in a set of RBM simulations, Elliott et al. [10] showed that head impacts are influenced by the timing of the vehicle impact relative to the pedestrian gait cycle, see Figure 11. The gait affects the upper body rotation, which ultimately affects the head impact velocities over the gait cycle, as well as the head impact location on the vehicle. The gait's influence on head impact response and impact location in lateral collisions was further confirmed in a study by Peng et al. [41], who studied a similar simulation setup.

### 4.5.2. Subject orientation angle

Several researchers have demonstrated the importance of the pedestrian's orientation angle relative to the vehicle for head injury prediction [9, 35, 91]. Tamura et al. [91] performed a series of pedestrian collisions using an HBM and reported that the pre-impact pedestrian orientation angle considerably affected the contact force, head velocity and CSDM. Coley et al. [9] reconstructed a real-world car-to-pedestrian accident and, in agreement with Tamura et al. [91], stated that the pedestrian orientation angle relative to the car at impact is one of the main parameters behind estimations of HIC values. Liu et al. [35], who performed a parameter sensitivity analysis of a simulated lateral pedestrian impact, also observed a strong association between estimated HIC values and the pedestrian orientation angle.

### 4.5.3. Subject initial position relative to vehicle

In this context, the subject initial position relative to the vehicle relates to the translational position of the pedestrian in front of the vehicle (how much right or left of the vehicle the pedestrian is located).

Coley et al. [9] reconstructed a real-world car-to-pedestrian accident using RBMs and urged that the pedestrian's initial position relative to the car prior to impact is one of the main parameters behind estimations of head injury risk. The subject position relative to the vehicle would dictate where on the vehicle the head would strike, which in turn will affect the force behind the impact (recall Figure 8).

### 4.5.4. Pre-crash reaction

Soni et al. [90] performed 40 volunteer experiments, where the volunteer was subjected to a simulated accident situation on a two-way traffic street-crossing. The volunteer's postural changes, speed and orientation was documented using Vicon markers, compiling a dataset of accident avoidance strategies, such as running, stopping and stepping back. Based on the dataset, Soni et al. then performed RBM simulations of car-to-pedestrian impacts, with the initial conditions based on the 40 obtained volunteer reactions. These simulations were compared with a baseline simulation, where the pedestrian model was positioned in a typical walking posture. For all simulations, the same car model was used, striking the pedestrian at 40 km/h. The results suggested that the avoidance reaction had an influence on the head impact location, head kinematics and HIC. However, the position of the struck-side leg (whether or not the left or right leg was leading prior to impact) was proven far more influential on the injury predictions than the crash avoidance reaction.

### 4.5.5. Muscle activation/tonus

Alvarez et al. [8] used FE simulation to investigate the effect of neck tonus on a pedestrian-to-vehicle impact on a generalized hood. Head impacts against a car was simulated using a pedestrian HBM with a relaxed neck, and subsequently compared with simulations in which the neck muscle tonus of the HBM was increased. It was concluded that the neck muscle tonus had a negligible effect on the predicted brain strains (1 to 14% change in MPS). The muscle tonus could however influence the head rotation during a vehicle collision, which in turn could alter the orientation of the head prior to impact [92].

## 4.6. Impact kinematics
### 4.6.1. Pedestrian speed

Selected studies suggest that the speed of the pedestrian prior to lateral vehicle collisions has minor effect on the head impact response [10, 35]. Elliott et al. [10] performed a set of RBM simulations of lateral car-to-pedestrian collisions, varying the vehicle speed, pedestrian speed and pedestrian gait. While the pedestrian speed had influence on the head impact location, the head rotation and head impact velocity was largely unaffected.

### 4.6.2. Vehicle speed

There is convincing amount of evidence suggesting that vehicle impact speed is one of the most influential covariates affecting head injury outcome in car-to-pedestrian collisions [5–7, 10, 17, 32, 34–36, 39, 41, 46, 51, 68, 71, 74, 93]. The impact speed reportedly dictates the magnitude of HIC values [5–7] and head impact velocities [93].

Looking at records of pedestrian accidents that occurred during the years 1999-2007, comprising 2127 car-to-pedestrian collisions, a strong dependence on impact speed and fatality risk was found. The risk of fatality at 50 km/h impacts was found to be more than twice the risk at 40 km/h, and more than five times higher than the risk at 30 km/h [94].

### 4.6.3. Vehicle steering/yaw angle

Researchers have pointed out that the steering angle is influential for the head-to-vehicle impact outcome. Wang et al. [7] used RBMs to investigate the influence of vehicle steering angle in car-to-pedestrian collisions and concluded that different vehicle steering and evasive maneuvers result in more lateral head collision forces and accelerations, which may influence head impact severity.

Pan et al. [95] performed a parametric study of vehicle steering maneuvers in simulations of a 40 km/h cyclist collision using whole-body HBMs. The simulation matrix included four vehicle yaw angular velocities and five different impact locations. It was proven how the yaw angle had a significant effect on the WAD, HIC, CSDM and head impact kinematics. The yaw angle effect on HBM kinematics and head angular accelerations is illustrated in Figure 12.

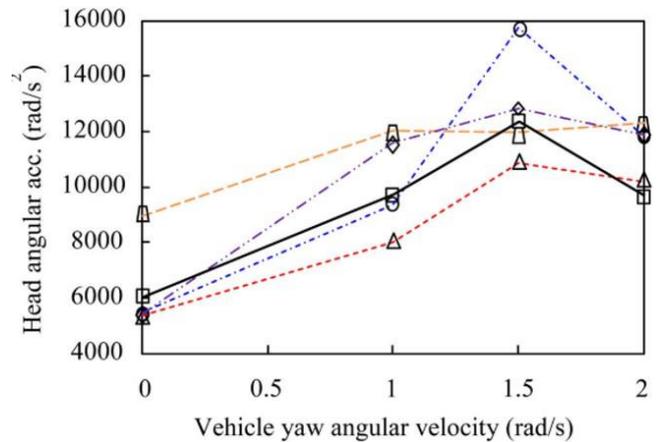

**Figure 12.** Illustrative figures showing how vehicle yaw angle may affect the body kinematics and head PAA in a cyclist impact. The different graphs illustrate different relative starting positions of the pedestrian. The figures are borrowed from Pan et al. [95] with permission from the authors.

### 4.6.4. Vehicle braking deceleration

In a study using headform impacts on a vehicle hood, Fredriksson et al. [71] was able to show how a 10 km/h reduction in impact speed, corresponding to the deceleration of an automatic braking system, could decrease the HIC

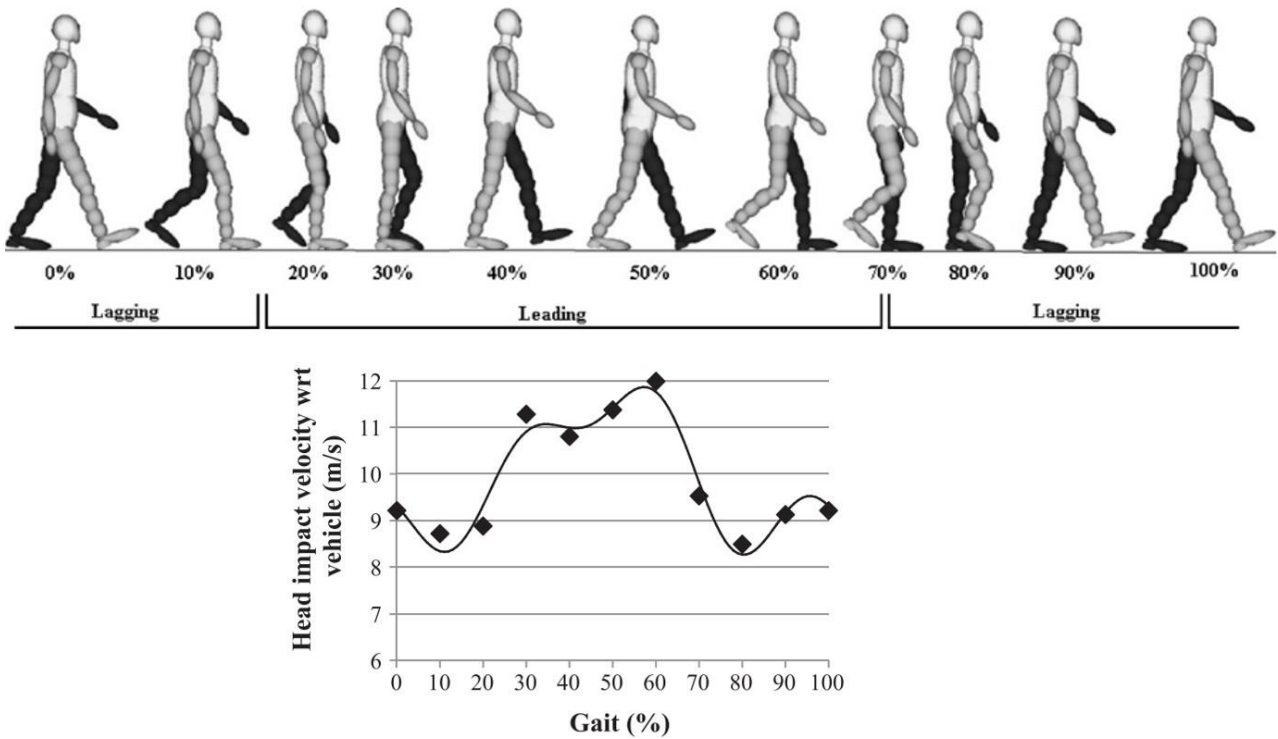

**Figure 11.** Top: The phases of the gait cycle. Bottom: Head impact velocity and rotation of the head at impact changes over the course of the gait cycle in lateral vehicle impacts. Figure borrowed from Elliott et al. [10] with permission from the authors.

values by up to 92%, followed by significant reductions in predicted brain strains.

Kawabe et al. [51] studied lateral impacts of a pedestrian HBM and showed that deceleration had an effect on the head trajectory and thus the head impact location (quantified by the WAD). The head kinematics can be affected as well. In Figure 13, the HBM's head CoG velocity curves with and without deceleration of an impacting vehicle (starting from 40 km/h) is presented.

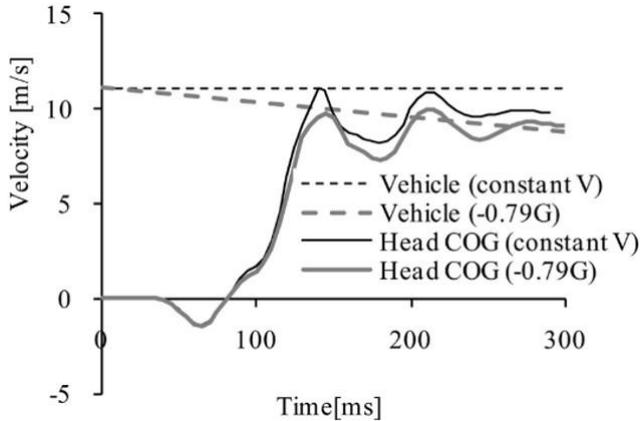

**Figure 13.** HBM head CoG velocity during a lateral car impact at 40 km/h, with and without an applied braking deceleration of the striking vehicle. The figure is borrowed from Kawabe et al. [51].

### 4.6.5. Vehicle pitching due to braking

Researchers have reported that pitching due to vehicle braking has slight effect on the head impact location or WAD, and insignificant effect on head kinematics [96, 97].

This conclusion was drawn by Fredriksson et al. [96], who modeled vehicle pitching by introducing a 1° rotation of the impacting vehicle in a study using FE simulation and pedestrian dummy experiments. By comparing the impact with and without the 1° pitch, it was shown how vehicle pitching due to braking has slight effect on the head impact location, but in the cases where the head impact locations were the same, no significant influence of the pitching was observed in head kinematics.

### 4.6.6. Pedestrian sliding

Several researchers point out that the head impact kinematics and impact point will be influenced by the pedestrian sliding/slipping up the hood [36, 55–58]. An illustration of pedestrian slippage/sliding is seen in Figure 14. Seemingly, it is not fully understood yet why and when sliding/slipping along the hood occurs.

Paas et al. [58] hypothesizes that the amount of sliding over the hood is determined by the pelvis-to-BLE height ratio. For large ratios, the subject's mass come into play as well. Based on PMHS tests, Paas et al. observed that subjects with a highly situated pelvis (as for tall subjects) and larger mass displayed more sliding compared to shorter and lighter subjects. The sliding had a pronounced effect on the head impact kinematics and impact locations. It was further claimed that sliding increased the head impact velocity by reducing the momentum transfer between the car and pedestrian. Kerrigan et al. [57] has hypothesized that the sliding distance could be further influenced by clothing, bonnet slope angle, BLE shape, friction and damage pattern on the BLE and hood. Liu et al. [36] hypothesizes that the sliding/slipping movement is connected with the knee lateral deflection. A large knee lateral deflection would "restrict" the slipping movement and exaggerate the rotational movement of the upper body towards the hood top.

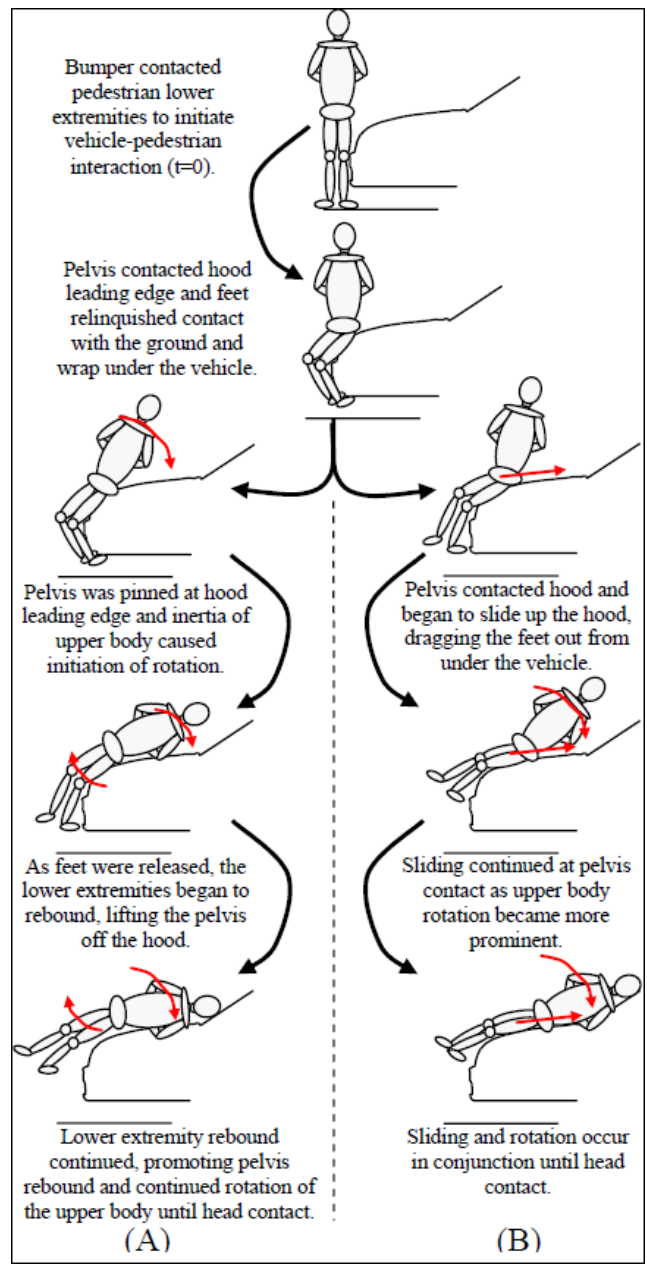

**Figure 14.** Illustration of sliding/slipping phenomenon in pedestrian lateral impacts: no slippage to the left, slippage to the right. Observe how the slipping changes the impact point of the head. Figure borrowed from Kerrigan et al. [55].

### 4.7. Boundary conditions
#### 4.7.1. Release time
Schroeder et al. [42] performed impact tests on four PMHS. One of the four subjects were released slightly earlier than the other subjects (the authors do not present the release time difference in quantitative terms). This resulted in different torso rotation compared to the other subjects, affecting the head trajectory and consequently, the head impact point on the car. Based on these findings, it can be hypothesized that the release time of the pedestrian is important, also in numerical simulations.

#### 4.7.2. Relative ground level
Raising or dropping the pedestrian ground level relative to the vehicle will directly influence the body kinematics and the head impact location of a pedestrian. The relative ground level will determine the location of the knee and hip joints with respect to the vehicle. Furthermore, it will determine the location of the pedestrian's CoG relative to the vehicle [87]. This is important in the context of this report since, among other things, the kinematics of the pelvis has been shown to be an important determinant of the overall upper body kinematics, including the head [55, 88].

#### 4.7.3. Car-to-pedestrian friction
Elliott et al. [50] studied the head's trajectory, translation, impact location, impact time and impact velocity in a series of RBM simulations of lateral car-to-pedestrian impacts. The researchers observed no changes in the head impact conditions when the friction coefficient between the vehicle and pedestrian was doubled or halved from a value of 0.3.

#### 4.7.4. Vehicle mass
Vehicles of different classes and models may have different masses. However, mass variation alone is not a particularly important covariate for car-to-pedestrian injury prediction. It has been shown that the mass difference between smaller cars and SUVs are negligible for pedestrian injury prediction, as the mass difference would result in almost no change in momentum transfer in the primary impact [6, 38, 43]. Meaning, the substantial mass difference between a passenger car of 1200 kg and a SUV of 1600 kg, and even larger trucks, will not have notable effect on pedestrian injury risk.

#### 4.7.5. Influence of other body parts
By comparing full-body HBMs with head-only FE models, it has been asserted that the body parts below the head, in particular the neck, are influential for the head impact kinematics [98, 99].

Elbow and shoulder impact have an effect on head rotation in lateral pedestrian collisions. Based on observations from four PMHS tests, Paas et al. [100] point out that elbows and shoulders provide additional support for the thorax and consequently change neck boundary conditions, increasing rotational accelerations towards the impacting vehicle.

Ishikawa et al. [101] compared RBM simulations with PMHS tests and suggested that the influence of the lower body motion to the upper body response seems insignificant during the initial car-to-lower-extremities contact, in terms of trajectories. The researchers further claim that, during the period between initial contact and the body rotation towards the hood, the upper extremities have an influence on the movement of the upper body.

Studying PMHS impacts, Paas et al. [58] observed that the upper arm and shoulder altered the head kinematics. The angle of the upper arm in lateral impacts had an influence on the head velocity.

Both Kerrigan et al. [55] and Subit et al. [88] concluded that the kinematics of the pelvis was an important determinant of the overall upper body kinematics. In turn, the pelvic kinematics response is determined by its height. These conclusions were drawn based on observations from PMHS impact experiments.

### 4.8. Knowledge gaps
Despite significant advances on this topic, several areas remain underexplored, offering opportunities for future research. The authors identified a number of potential covariates that were not brought up in the body of literature reviewed in this study. These knowledge gaps are presented in this section and are also highlighted in Figure 4, written in italic font.

#### 4.8.1. Body composition and HBM mass distribution
Usually, when an HBM is being scaled or morphed, the HBM is globally reshaped, without deforming fat tissue and muscle tissue, and sometimes even bone, separately. Many authors also tune the density of the soft tissues uniformly of an HBM to adjust its total mass. After morphing, scaling or density tuning, no attention is usually given to the resulting CoG or mass distribution or the HBM. In reality, the same BMI pedestrian can have very different body compositions and body mass distributions (Heymsfield et al. 2007). It has also been reported that tissue properties can vary with age and sex. For example, bone density will generally decrease over time, and age-related bone loss can be more pronounced in females compared to males (Mumtaz et al., 2020). Ultimately, whether the HBM mass distribution (CoG), which can be influenced by variations in body composition, sex and age, has a non-negligible contribution to pedestrian impact response has not been covered in literature. No authors have so far covered or proven whether or not the HBM mass distribution, or the HBM's CoG, would explicitly influence the predicted head impact response.

#### 4.8.2. Choice of statistical body shapes
The 50$^{th}$ percentile American/European male is the most common target anthropometry of some of the most frequently used HBMs, such as THUMS, GHBMC and SAFER. In the processed literature, most authors choose an American/European 50$^{th}$ percentile male for their analysis. If geometrical morphing is performed to achieve a targeted

BMI, statistical body shapes based on American/European bodies have mainly been used, e.g. using UMTRI [102] or SMPL [103]. It is not well-known if statistical body shapes of equal BMI differ significantly among different datasets, and whether or not the choice of statistical body model subsequently would have a significant influence on injury predictions in car-to-pedestrian impacts.

### 4.8.3. Choice of FE models

Researchers in the field all use their choice of HBM, and there are nowadays several to choose from (THUMS, GHBMC, VIVA+, SAFER etc). Whether the choice of HBM is an important covariate in predicting head injuries in car-to-pedestrian impacts remains unclear. Likewise, there are many available FE vehicle models being used in FE accident reconstructions, some of them simplified to modular car models, while some of them are as detailed as a real-world vehicle. Assuming the front-end profile of the vehicles are similar, would the choice of vehicle model have a significant influence on the head impact predictions? An area that remains underexplored is whether the choice of HBM or vehicle model have a significant effect on head injury predictions.

### 4.8.4. Settling of HBMs and feet-to-ground friction

It seems to be praxis to settle an HBM to the ground before simulating an impact. The settling implicates allowing your HBM to rest under gravity for a short time window prior to impact, settling the feet on the ground. No authors seem to cover whether or not this simulation step is actually important. Is the ground and the HBM settling necessary for accurate head injury predictions of head-to-vehicle impacts? This goes hand-in-hand with the question of whether the feet-to-ground friction is influential to head impacts, a topic that neither has been addressed in literature.

### 4.8.5. Ground clearance

Researchers have investigated the influence of several car-related measurements, including the height of the BLE, the length of the hood, the height and lead of the bumper height and so forth. However, to the authors knowledge, no researcher has investigated the influence of the ground clearance, meaning the distance between the ground and the lowest point of the vehicle front. Hypothetically, the risk with a too large ground clearance is that the feet and ankles might wrap under the vehicle, hindering any sliding/slipping on the vehicle hood, or change the rotation around the pelvis. Investigating whether or not the ground clearance should be modeled with care or not in FE accident reconstructions of vehicle-to-pedestrian collisions could be relevant.

### 4.8.6. Windshield boundary conditions

It has been shown that the windshield frame/borders are notably stiff and particularly hazardous to pedestrians, see Figure 8. However, no authors scrutinize how this windshield frame should be modeled in FE. Are rigid boundary conditions for the windshield edge, which is suspected to be used by most authors, valid? The boundary conditions of the windshield might influence the initiation and pattern of predicted windshield cracks, and should thus be investigated.

## 5. Discussion

In this review, variables reported to have an influence on head injury predictions in vehicle-to-pedestrian collisions have been identified and listed. The listed variables, referred-to as covariates, relate to the modeling of both the impacting vehicle and the impacted pedestrian, as well as the impact conditions. The summary of covariates presented in this review can be used to make more informative decisions when carrying out VRU accident reconstructions with FE techniques, helping researchers to avoid making simplifications of FE models that compromise the accuracy of numerical head injury predictions. The presented material could also be used when selecting parameters to include in sensitivity analyses of car-to-pedestrian impact simulations.

The summary of findings presented in this review highlight important aspects of conducting car-to-pedestrian accident reconstructions. First of all, it is imperative to take a subject's anthropometry and body posture into account- generic, unpositioned HBMs are simply not adequate to capture pedestrian head traumas correctly. Second, head injuries in a car-to-pedestrian accident can simply not be accurately predicted using generic vehicle models- the car's front-end geometry, its stiffness and construction plays an outsized role in predictions of head impact response. Thirdly, the impact conditions must be thoroughly derived- the impact velocity, braking and yaw angle is far too influential to be roughly estimated. On the other hand, this review also showed, that a modeler aiming to predict pedestrian head impact response using FE can pay less attention to other aspects, such as precisely modeling the bumper of an impacting vehicle.

A key limitation of this study is the relatively small set of literature included in the review. Some relevant research may have been overlooked, and the conclusions drawn may not fully represent the breadth of the field. However, despite this limitation, the review highlights some of the most important findings from the relevant studies available. The purpose of the study was to provide an informative overview of the current state of research on covariates in pedestrian accidents rather than being all-encompassing. While not exhaustive, the selected literature captures the core insights and trends in the field, offering valuable contributions to the understanding of the topic.

It should also be mentioned, that the conclusions drawn by the studies referenced in this review were primarily drawn from studies using RBMs rather than HBMs. While RBM simulations can provide useful insights, they fail to capture the complex, deformable nature of the human body, limiting their accuracy in representing real-world biomechanics. In consequence, many of the references studies only evaluate kinematics-based injury risk metrics, such as HIC, while a small subset of them also evaluate tissue-based metrics, such as MPS. Therefore, whether some of the listed covariates actually influence brain deformation has not been fully proven.

This review has highlighted several important knowledge gaps within the current literature, which present valuable opportunities for advancing the field. The listed gaps remain underexplored or inconsistently addressed. By focusing future research on these areas, researchers can help refine existing models and improve practical applications. Addressing these gaps will not only enhance the knowledge but also lead to more accurate and efficient accident reconstructions. The identification of these gaps serves as a roadmap for future studies, providing a direction for advancing research in this area.

## 6. Acknowledgments

This study was partly financed by the Swedish Research Council (VR, Nos. 2020-04724 and 2020-04496) and KTH Digital Futures Research Pairs Project (2024-2025).

## References

[1] M. C. Dewan, A. Rattani, S. Gupta, R. E. Baticulon, Y. C. Hung, M. Punchak, A. Agrawal, A. O. Adeleye, M. G. Shrime, A. M. Rubiano, J. V. Rosenfeld, K. B. Park, Estimating the global incidence of traumatic brain injury, Journal of Neurosurgery 130 (2019) 1080–1097.

[2] A. I. Maas, D. K. Menon, G. T. Manley, M. Abrams, C. Åkerlund, N. Andelic, et al., Traumatic brain injury: progress and challenges in prevention, clinical care, and research, The Lancet Neurology 21 (2022) 1004–1060.

[3] World Health Organization, Global status report on road safety 2023, Technical Report, 2023. URL: https://iris.who.int/bitstream/handle/10665/375016/9789240086517-eng.pdf?sequence=1.

[4] J. R. Crandall, D. Bose, J. Forman, C. D. Untaroiu, C. Arregui-Dalmases, C. G. Shaw, J. R. Kerrigan, Human surrogates for injury biomechanics research, Clinical Anatomy 24 (2011) 362–371.

[5] R. Kendall, M. Meissner, J. Crandall, The causes of head injury in vehicle-pedestrian impacts: Comparing the relative danger of vehicle and road surface, SAE Technical Papers (2006).

[6] C. K. Simms, D. P. Wood, Pedestrian risk from cars and sport utility vehicles - A comparative analytical study, Proceedings of the Institution of Mechanical Engineers, Part D: Journal of Automobile Engineering 220 (2006) 1085–1100.

[7] D. Wang, W. Deng, L. Wu, L. Xin, L. Xie, H. Zhang, Impact of Vehicle Steering Strategy on the Severity of Pedestrian Head Injury (2024).

[8] V. S. Alvarez, P. Halldin, S. Kleiven, The Influence of Neck Muscle Tonus and Posture on Brain Tissue Strain in Pedestrian Head Impacts, SAE Technical Papers 2014-Novem (2014) 1–40.

[9] G. Coley, R. D. Lange, P. D. Oliveira, R. Happee, Pedestrian Human Body Validation Using Detailed Real-World Accidents ircobi, Ircobi (2001).

[10] J. R. Elliott, C. K. Simms, D. P. Wood, Pedestrian head translation, rotation and impact velocity: The influence of vehicle speed, pedestrian speed and pedestrian gait, Accident Analysis and Prevention 45 (2012) 342–353.

[11] J. Wang, Z. Li, F. ying, D. Zou, Y. Chen, Reconstruction of a real-world car-to-pedestrian collision using geomatics techniques and numerical simulations, Journal of Forensic and Legal Medicine 91 (2022) 102433.

[12] M. Lalwala, A. Chawla, P. Thomas, S. Mukherjee, Finite element reconstruction of real-world pedestrian accidents using THUMS pedestrian model, International Journal of Crashworthiness 25 (2020) 360–375.

[13] N. Trube, P. Matt, M. Jenerowicz, N. Ballal, T. Soot, D. Fressmann, N. Lazarov, J. Moennich, T. Lich, P. Lerge, L. V. Nölle, S. Schmitt, Plausibility Assessment of Numerical Cyclist to Vehicle Collision Simulations based on Accident Data, in: Proceedings of the annual International Research Council on the Biomechanics of Injury (IRCOBI), Cambridge, UK, 2023, pp. 113–135.

[14] C. Costa, J. Aira, B. Koya, W. Decker, J. Sink, S. Withers, R. Beal, S. Schieffer, S. Gayzik, J. Stitzel, A. Weaver, Finite element reconstruction of a vehicle-to-pedestrian impact, Traffic Injury Prevention 21 (2020) S145–S147.

[15] B. Benea, A. Soica, The Contact Phase in Vehicle–Pedestrian Accident Reconstruction, Applied Sciences (Switzerland) 13 (2023).

[16] A. Tamura, T. Koide, K. H. Yang, Effects of ground impact on traumatic brain injury in a fender vault pedestrian crash, International Journal of Vehicle Safety 8 (2015) 85–100.

[17] R. Fredriksson, E. Rosén, A. Kullgren, Priorities of pedestrian protection - A real-life study of severe injuries and car sources, Accident Analysis and Prevention 42 (2010) 1672–1681.

[18] B. S. Roudsari, C. N. Mock, R. Kaufman, An evaluation of the association between vehicle type and the source and severity of pedestrian injuries, Traffic Injury Prevention 6 (2005) 185–192.

[19] M. J. Aminoff, F. Boller, D. F. Swaab, Handbook of Clinical Neurology, volume 127, Elsevier B.V, Amsterdam, Netherlands, 2015.

[20] F. A. Fernandes, R. J. D. Sousa, Head injury predictors in sports trauma - A state-of-the-art review, Proceedings of the Institution of Mechanical Engineers, Part H: Journal of Engineering in Medicine 229 (2015) 592–608.

[21] J. A. Newman, Head injury criteria in automotive crash testing, SAE Technical Papers 89 (1980) 701–747.

[22] S. Kleiven, Why Most Traumatic Brain Injuries are Not Caused by Linear Acceleration but Skull Fractures are, Frontiers in Bioengineering and Biotechnology 1 (2013) 1–5.

[23] J. A. Newman, A Generalized Model for Brain Injury Threshold, In Proceedings of International Conference on the Biomechanics of Impact, 1986 (pp. 121-131). (1986) 121–131.

[24] J. A. Newman, N. Shewchenko, A Proposed New Biomechanical Head Injury Assessment Function - The Maximum Power Index, SAE Technical Papers 2000-Novem (2000).

[25] E. G. Takhounts, M. J. Craig, K. Moorhouse, J. McFadden, V. Hasija, Development of Brain Injury Criteria (BrIC), SAE Technical Papers 2013-Novem (2013) 243–266.

[26] L. F. Gabler, J. R. Crandall, M. B. Panzer, Development of a Metric for Predicting Brain Strain Responses Using Head Kinematics, Annals of Biomedical Engineering 46 (2018) 972–985.

[27] J. R. Crandall, K. S. Bhalla, N. J. Madeley, Designing road vehicles for pedestrian protection, British Medical Journal 324 (2002) 1145–1148.

[28] C. Simms, D. Wood, Pedestrian and Cyclist Impact: A Biomechanical Perspective, volume 55, Springer Science+Business Media, B.V., 2010. doi:10.1111/j.1556-4029.2010.01360.x.

[29] R. W. Anderson, S. Doecke, An analysis of head impact severity in simulations of collisions between pedestrians and SUVs/Work utility vehicles, and sedans, Traffic Injury Prevention 12 (2011) 388–397.

[30] A. J. Fisher, R. R. Hall, Accident details Accident sample Pedestrian Car make Number of models, Accident Analysis and Prevention 4 (1972) 47–58.

[31] V. Gupta, K. H. Yang, Effect of Vehicle Front End Profiles Leading to Pedestrian Secondary Head Impact to Ground, SAE Technical Papers 2013-Novem (2013) 139–155.

[32] Y. Han, J. Yang, K. Mizuno, Y. Matsui, Effects of Vehicle Impact Velocity, Vehicle Front-End Shapes on Pedestrian Injury Risk, Traffic Injury Prevention 13 (2012) 507–518.

[33] Y. Han, J. Yang, K. Nishimoto, K. Mizuno, Y. Matsui, D. Nakane, S. Wanami, M. Hitosugi, Finite element analysis of kinematic behaviour and injuries to pedestrians in vehicle collisions, International Journal of Crashworthiness 17 (2012) 141–152.

[34] H. Li, K. Li, W. Lv, S. Cui, L. He, S. Ruan, Analysis of the effects of vehicle model and speed on the head injury of six-year-old pedestrian, International Journal of Vehicle Safety 12 (2021) 177–193.

[35] W. Liu, A. Duan, K. Li, J. Qiu, L. Fu, H. Jia, Z. Yin, Parameter sensitivity analysis of pedestrian head dynamic response and injuries based on coupling simulations, Science Progress 103 (2020) 1–15.

[36] X. J. Liu, J. K. Yang, P. Lövsund, A study of influences of vehicle speed and front structure on pedestrian impact responses using mathematical models, Traffic Injury Prevention 3 (2002) 31–42.

[37] D. Longhitano, B. Henary, K. Bhalla, J. Ivarsson, J. Crandall, Influence of vehicle body type on pedestrian injury distribution, SAE Technical Papers (2005).

[38] K. Mizuno, J. Kajzer, Compatibility problems in frontal, side, single car collisions and car-to-pedestrian accidents in Japan, Accident Analysis and Prevention 31 (1999) 381–391.

[39] K. Mizuno, J. Kajzer, Head Injuries in Vehicle-Pedestrian Impact, SAE International (2000).

[40] Y. Okamoto, T. Sugimoto, K. Enomoto, J. Kikuchi, Pedestrian head impact conditions depending on the vehicle front shape and its construction - Full model simulation, Traffic Injury Prevention 4 (2003) 74–82.

[41] Y. Peng, C. Deck, J. Yang, R. Willinger, Effects of pedestrian gait, vehicle-front geometry and impact velocity on kinematics of adult and child pedestrian head, International Journal of Crashworthiness 17 (2012) 553–561.

[42] G. Schroeder, K. Fukuyama, K. Yamazaki, K. Kamiji, T. Yasuki, Injury mechanism of pedestrians impact test with a sport-utility vehicle and mini-van, International Research Council on the Biomechanics of Injury - 2008 International IRCOBI Conference on the Biomechanics of Injury, Proceedings (2008) 259–273.

[43] C. K. Simms, D. Wood, R. Fredriksson, Accidental Injury: Biomechanics and Prevention, in: Accidental Injury: Biomechanics and Prevention, Springer Science+Business Media New, New York, 2015, pp. 721–753. doi:10.1007/978-1-4939-1732-7.

[44] A. Tamura, K. H. Yang, Essential factors leading to a traumatic brain injury during low-speed fender vault pedestrian impacts, International Journal of Vehicle Safety 11 (2019) 1–18.

[45] B. A. Tolea, H. Beles, A. I. Radu, F. Scurt, G. Dragomir, The assessment of pedestrian's head injury risk at the impact with the vehicle's windshield, IOP Conference Series: Materials Science and Engineering 1256 (2022) 012042.

[46] J. Yao, J. Yang, D. Otte, Head injuries in child pedestrian accidents - In-depth case analysis and reconstructions, Traffic Injury Prevention 8 (2007) 94–100.

[47] G. Zhang, Q. Qin, Z. Chen, Z. Bai, L. Cao, A study of the effect of the front-end styling of sport utility vehicles on pedestrian head injuries, Applied Bionics and Biomechanics 2018 (2018).

[48] J. Kerrigan, C. Arregui-Dalmases, J. Crandall, Assessment of pedestrian head impact dynamics in small sedan and large SUV collisions, International Journal of Crashworthiness 17 (2012) 243–258.

[49] L. Martinez, L. J. Guerra, G. Ferichola, A. Garcia, J. Yang, J. Yao, Stiffness corridors for the current European Fleet, APROSYS SP3 Deliverable Report D312B. (2006).

[50] J. R. Elliott, M. Lyons, J. Kerrigan, D. P. Wood, C. K. Simms, Predictive capabilities of the MADYMO multibody pedestrian model: Three-dimensional head translation and rotation, head impact time and head impact velocity, Proceedings of the Institution of Mechanical Engineers, Part K: Journal of Multi-body Dynamics 226 (2012) 266–277.

[51] Y. Kawabe, T. Asai, D. Murakami, C. Pal, T. Okabe, Different Factors Influencing Post-crash Pedestrian Kinematics, SAE International Journal of Passenger Cars - Mechanical Systems 5 (2012) 214–230.

[52] L. Shi, Y. Han, H. Huang, W. He, F. Wang, B. Wang, Effects of vehicle front-end safety countermeasures on pedestrian head injury risk during ground impact, Proceedings of the Institution of Mechanical Engineers, Part D: Journal of Automobile Engineering 233 (2019) 3588–3599.

[53] J. Yang, P. Lövsund, C. Cavallero, J. Bonnoit, A Human-Body 3D Mathematical Model for Simulation of Car-Pedestrian Impacts, Journal of Crash Prevention and Injury Control 2 (2000) 131–149.

[54] S. Yin, J. Li, J. Xu, Exploring the mechanisms of vehicle front-end shape on pedestrian head injuries caused by ground impact, Accident Analysis and Prevention 106 (2017) 285–296.

[55] J. R. Kerrigan, J. R. Crandall, B. Deng, Pedestrian kinematic response to mid-sized vehicle impact, International Journal of Vehicle Safety 2 (2007) 221–240.

[56] R. W. G. Anderson, L. D. Streeter, G. Ponte, J. McLean, Mc, Pedestrian Reconstruction Using Multibody Madymo Simulation And The Polar-li Dummy: A Comparison Of Head Kinematics, ESV -Paper no. 07-0273 (2007) 1–15.

[57] J. R. Kerrigan, D. B. Murphy, D. C. Drinkwater, C. Y. Kam, D. Bose, J. R. Crandall, Kinematic corridors for PMHS tested in full-scale pedestrian impact tests, in: ESV Confer- ence, Washington, DC, 2005, 2014, pp. 1–18.

[58] R. Paas, C. Masson, J. Davidsson, Head boundary conditions in pedestrian crashes with passenger cars: Six-degrees-of-freedom post-mortem human subject responses, International Journal of Crashworthiness 20 (2015) 547–559.

[59] G. Stcherbatcheff, C. Tarriere, P. Duclos, A. Fayon, Simulation of Collisions Between Pedestrians and Vehicles Using Adult and Child Dummies, SAE Technical Papers (1975).

[60] A. Ahmed, The influence of the vehicle hood inclination angle on the severity of the pedestrian adult head injury in a front collision using finite element modeling, Thin-Walled Structures 150 (2020).

[61] B. Tolea, A. I. Radu, H. Beles, C. Antonya, Influence of the geometric parameters of the vehicle frontal profile on the pedestrian's head accelerations in case of accidents, International Journal of Automotive Techonology 19 (2018) 85–98.

[62] Z. Cai, Y. Xia, X. Huang, ANALYSES of PEDESTRIAN'S HEAD-TO-WINDSHIELD IMPACT BIOMECHANICAL RESPONSES and HEAD INJURIES USING A HEAD FINITE ELEMENT MODEL, Journal of Mechanics in Medicine and Biology 20 (2020) 1–15.

[63] G. Li, M. Lyons, B. Wang, J. Yang, D. Otte, C. Simms, The influence of passenger car front shape on pedestrian injury risk observed from German in-depth accident data, Accident Analysis and Prevention 101 (2017) 11–21.

[64] F. Wang, M. Wang, L. Hu, K. Peng, J. Yin, D. Wang, L. Shi, Z. Zhou, Effects of the windshield inclination angle on head/brain injuries in car-to-pedestrian collisions using computational biomechanics models, Transportation Safety and Environment 6 (2024).

[65] M. Lyons, C. K. Simms, Predicting the influence of windscreen design on pedestrian head injuries, 2012 IRCOBI Conference Proceedings - International Research Council on the Biomechanics of Injury (2012) 703–716.

[66] A. Konosu, H. Ishikawa, A. Sasaki, A study on pedestrian impact procedure by computer simulation, Sustainability (Switzerland) 11 (1998) 1–14.

[67] A. Otubushin, D. Hughes, N. Ridley, K. Taylor, Analysis of pedestrian head impacts to the bonnets of European vehicles, International Journal of Crashworthiness 4 (1999) 159–174.

[68] R. Watanabe, T. Katsuhara, H. Miyazaki, Y. Kitagawa, T. Yasuki, Research of the Relationship of Pedestrian Injury to Collision Speed, Car-type, Impact Location and Pedestrian Sizes using Human FE model (THUMS Version 4), SAE Technical Papers 2012-Octob (2012).

[69] X. Lv, W. Liu, D. Zhou, C. Wang, F. Ma, F. Zhao, Assessment and optimization of pedestrian head impact safety performance for a SUV car, Proceedings - 2013 5th Conference on Measuring

Technology and Mechatronics Automation, ICMTMA 2013 (2013) 313–316.
[70] K. Mizuno, H. Yonezawa, J. Kajzer, Pedestrian headform impact tests for various vehicle locations, Proceedings of the 17th Enhanced Safety Vehicle (ESV) Conference (2001) 1–10.
[71] R. Fredriksson, L. Zhang, O. Boström, K. Yang, Influence of Impact Speed on Head and Brain Injury Outcome in Vulnerable Road User Impacts to the Car Hood, SAE Technical Papers 2007-Octob (2007) 155–167.
[72] R. Fredriksson, L. Zhang, O. Böstrom, Influence of deployable hood systems on finite element modelled brain response for vulnerable road users, International Journal of Vehicle Safety 4 (2009) 29–44.
[73] R. Fredriksson, Y. Håland, J. Yang, Evaluation of a new pedestrian head injury protection system with a sensor in the bumper and lifting of the bonnet's rear part, Lecture Notes in Engineering and Computer Science 3 LNECS (2013) 1947–1952.
[74] N. Iwai, T. Araki, Numerical simulation of pedestrian head impact on vehicle front structure, SAE Technical Papers 3 (2003) 653–658.
[75] T. F. MacLaughlin, J. W. Kessler, Pedestrian head impact against the central hood of motor vehicles - Test procedure and results, SAE Technical Papers 99 (1990) 1729–1737.
[76] H. B. Pritz, Experimental Investigation of Pedestrian Head Impacts on Hoods and Fenders of Production Vehicles (1983).
[77] H. Chen, J. Crandall, M. Panzer, Evaluating pedestrian head Sub-System test procedure against full-scale vehicle-pedestrian impact, International Journal of Crashworthiness 26 (2021) 467–489.
[78] V. Gupta, K. H. Yang, Effect of hood periphery shape parameters and hood-fender interface characteristics on pedestrian head injury criterion, International Journal of Vehicle Safety 8 (2015) 233–250.
[79] M. Ptak, D. Czerwińska, J. Wilhelm, F. A. Fernandes, R. J. de Sousa, Head-to-bonnet impact using finite element head model, Lecture Notes in Mechanical Engineering 1 (2019) 545–555.
[80] B. Nie, Q. Zhou, Y. Xia, J. Tang, Influence of feature lines of vehicle hood styling on headform kinematics and injury evaluation in car-to-pedestrian impact simulations, SAE International Journal of Transportation Safety 2 (2014) 182–189.
[81] B. M. Boggess, J. Wong, S. Mark, Development of Plastic Components for Pedestrian Head Injury Risk Reduction, International Technical Conference on Enhanced Safety of Vehicles (2003) 1–8.
[82] N. Jayanth, A. Agarwal, C. Sekhar, A Methodology to Enhance the Directional Load Bearing Performance of Cowl Cover and Its Effect on Pedestrian Head Impact, SAE Technical Papers 2020-April (2020) 1–8.
[83] M. Munsch, Lateral Glazing Characterization Under Head Impact : (2007) 1–13.
[84] V. S. Alvarez, S. Kleiven, Importance of windscreen modelling approach for head injury prediction, in: Proceedings of the annual International Research Council on the Biomechanics of Injury (IRCOBI), Malaga, Spain, 2016, pp. 813–830.
[85] H. Chen, D. Poulard, J. Forman, J. Crandall, M. B. Panzer, Evaluation of geometrically personalized THUMS pedestrian model response against sedan–pedestrian PMHS impact test data, Traffic Injury Prevention 19 (2018) 542–548.
[86] R. Paas, J. Östh, J. Davidsson, Which pragmatic finite element human body model scaling technique can most accurately predict head impact conditions in pedestrian-car crashes?, 2015 IRCOBI Conference Proceedings - International Research Council on the Biomechanics of Injury (2015) 546–576.
[87] J. Shin, C. Untaroiu, J. Kerrigan, J. Crandall, D. Subit, Y. Takahashi, A. Akiyama, Y. Kikuchi, D. Longhitano, Investigating pedestrian kinematics with the polar-II finite element model, SAE Technical Papers 116 (2007) 655–667.
[88] D. Subit, J. Kerrigan, J. Crandall, K. Fukuyama, K. Yamazaki, K. Kamiji, T. Yasuki, Pedestrian-vehicle interaction: Kinematics and injury analysis of four full-scale tests, International Research Council on the Biomechanics of Injury - 2008 International IRCOBI Conference on the Biomechanics of Injury, Proceedings (2008) 275–294.
[89] L. Yan, C. Liu, X. Zhu, D. Zhou, X. Lv, X. Kuang, Translational medical bioengineering research of traumatic brain injury among Chinese and American pedestrians caused by vehicle collision based on human body finite element modeling, Frontiers in Neurology 14 (2024) 1–11.
[90] A. Soni, T. Robert, P. Beillas, IRC-13-92 IRCOBI Conference 2013 (2013) 762–776.
[91] A. Tamura, Y. Nakahira, M. Iwamoto, K. Nagayama, T. Matsumoto, Effects of pre-impact body orientation on traumatic brain injury in a vehicle – pedestrian collision Atsutaka Tamura * Yuko Nakahira and Masami Iwamoto Kazuaki Nagayama and Takeo Matsumoto 3 (2008).
[92] V. S. Alvarez, M. Fahlstedt, P. Halldin, S. Kleiven, Importance of neck muscle tonus in head kinematics during pedestrian accidents, 2013 IRCOBI Conference Proceedings - International Research Council on the Biomechanics of Injury (2013) 747–761.
[93] Y. Peng, C. Deck, J. Yang, D. Otte, R. Willinger, A Study of Adult Pedestrian Head Impact Conditions and Injury Risks in Passenger Car Collisions Based on Real-World Accident Data, Traffic Injury Prevention 14 (2013) 639–646.
[94] E. Rosén, U. Sander, Pedestrian fatality risk as a function of car impact speed, Accident Analysis and Prevention 41 (2009) 536–542.
[95] D. Pan, Y. Han, W. He, H. Huang, Effect of vehicle steering maneuvers on kinematics and head injury risks of cyclists via finite element modeling analysis, International Journal of Crashworthiness 26 (2021) 608–616.
[96] R. Fredriksson, J. Shin, C. D. Untaroiu, Potential of pedestrian protection systems-A parameter study using finite element models of pedestrian dummy and generic passenger vehicles, Traffic Injury Prevention 12 (2011) 398–411.
[97] A. Konosu, Reconstruction analysis for car-pedestrian accidents using a computer simulation model, JSAE Review 23 (2002) 357–363.
[98] M. Fahlstedt, P. Halldin, V. S. Alvarez, S. Kleiven, Influence of the Body and Neck on Head Kinematics and Brain Injury Risk in Bicycle Accident Situations, 2016 IRCOBI Conference Proceedings - International Research Council on the Biomechanics of Injury (2016) 459–478.
[99] F. Wang, C. Yu, B. Wang, G. Li, K. Miller, A. Wittek, F. Wang, C. Yu, B. Wang, G. Li, K. Miller, Prediction of pedestrian brain injury due to vehicle impact using computational biomechanics models: Are head-only models sufficient?, Traffic Injury Prevention 21 (2020) 102–107.
[100] R. Paas, J. Davidsson, C. Masson, U. Sander, K. Brolin, J. Yang, Pedestrian shoulder and spine kinematics in full-scale PMHS tests for human body model evaluation, 2012 IRCOBI Conference Proceedings - International Research Council on the Biomechanics of Injury (2012) 730–750.
[101] H. Ishikawa, J. Kajzer, Computer simulation of impact response of the human knee joint in car-pedestrian accidents, SAE Technical Papers (1992).
[102] N. Lindgren, M. J. Henningsen, C. Jacobsen, C. Villa, S. Kleiven, X. Li, Prediction of Real-Life Skull Fracture Patterns using Subject-Specific FE Head Models, in: CMBBE 2023 Symposium, May, Paris, France, 2023.
[103] N. Lindgren, Q. Huang, Q. Yuan, M. Lin, P. Wang, S. Kleiven, X. Li, Toward systematic finite element reconstructions of accidents involving vulnerable road users, Traffic Injury Prevention 0 (2025) 1–12.